\documentclass[11pt]{article}

\usepackage{amscd,amsmath, amssymb, fancyhdr, mathbbol}

\numberwithin{equation}{section}

\def\eqref#1{(\ref{#1})}

\newcommand{\arrow}{{\:\longrightarrow\:}}
\newcommand{\Z}{{\Bbb Z}}
\def\C{{\Bbb C}}
\newcommand{\R}{{\Bbb R}}

\def\1{\sqrt{-1}\:}

\newcommand{\cntrct}                
{\hspace{2pt}\raisebox{1pt}{\text{$\lrcorner$}}\hspace{2pt}}

\renewcommand{\bar}{\overline}
\renewcommand{\phi}{\varphi}
\renewcommand{\epsilon}{\varepsilon}
\renewcommand{\geq}{\geqslant}
\renewcommand{\leq}{\leqslant}

\newcommand{\Sp}{\operatorname{Sp}}
\newcommand{\U}{\operatorname{U}}

\newcounter{Mycounter}[section]
\newcounter{lemma}[section]
\newcounter{claim}[section]
\newcounter{sublemma}[section]
\newcounter{corollary}[section]
\newcounter{theorem}[section]
\renewcommand{\thetheorem}{{Theorem \thesection.\arabic{theorem}}}
\newcommand{\theorem}{%
    \setcounter{theorem}{\value{Mycounter}}
    \refstepcounter{theorem}
    \stepcounter{Mycounter}
    {\noindent \bf \thetheorem:\ }}

\newcounter{conjecture}[section]
\renewcommand{\theconjecture}{{Conjecture \thesection.\arabic{conjecture}}}
\newcommand{\conjecture}{%
    \setcounter{conjecture}{\value{Mycounter}}
    \refstepcounter{conjecture}
    \stepcounter{Mycounter}
    {\noindent \bf \theconjecture:\ }}

\newcounter{proposition}[section]
\renewcommand{\theproposition}
      {{Proposition \thesection.\arabic{proposition}}}
\newcommand{\proposition}{%
    \setcounter{proposition}{\value{Mycounter}}
    \refstepcounter{proposition}
    \stepcounter{Mycounter}
    {\noindent \bf \theproposition:\ }}

\newcounter{definition}[section]
\renewcommand{\thedefinition}
      {{Definition~\thesection.\arabic{definition}}}
\newcommand{\definition}{%
    \setcounter{definition}{\value{Mycounter}}
    \refstepcounter{definition}
    \stepcounter{Mycounter}
    {\noindent \bf \thedefinition:\ }}

\newcounter{example}[section]
\newcounter{remark}[section]
\renewcommand{\theremark}{{Remark \thesection.\arabic{remark}}}
\newcommand{\remark}{%
    \setcounter{remark}{\value{Mycounter}}
    \refstepcounter{remark}
    \stepcounter{Mycounter}
    {\noindent \bf \theremark:\ }}

\newcounter{problem}[section]
\newcounter{question}[section]
\makeatletter

\@addtoreset{equation}{section} \@addtoreset{footnote}{section}
\makeatother

\def\blacksquare{\hbox{\vrule width 5pt height 5pt depth 0pt}}
\def\endproof{\blacksquare}

\begin{document}
\begin{center}
{\LARGE\bf
Finiteness of stable Lagrangian fibrations
\\[4mm]
}

Ljudmila Kamenova
\footnote{This work was partially supported by a grant from the Simons 
Foundation/SFARI (522730, LK) 2017-2024}

\end{center}

\begin{center}{\em \small
In memory of my dear friend Sasha Ananin}
\end{center}

{\small \hspace{0.1\linewidth}
\begin{minipage}[t]{0.8\linewidth}
{\bf Abstract} \\
In this paper we survey some finiteness results of the deformation classes 
of hyperk\"ahler Lagrangian fibrations, and we prove finiteness 
for stable Lagrangian fibrations with a given discriminant divisor. 
\end{minipage}
}

{\scriptsize
\tableofcontents
}


\section{Introduction}


In hyperk\"ahler geometry a natural question to ask is whether there are 
finitely many examples (up to deformation) in a given dimension. For the 
known compact hyperk\"ahler examples, this is indeed the case. In each 
possible dimension, there are the classical examples of Hilbert schemes of 
points on a $K3$ surface, and the generalized Kummer varieties. In complex 
dimensions $6$ and $10$, there are additionally O'Grady's exceptional examples. 

\hfill

There are several finiteness results in the literature, establishing finiteness 
of deformations of various classes of hyperk\"ahler manifolds. Let us review 
some of them. 

\hfill

\theorem (Huybrechts, \cite{_Huybrechts:finiteness_}) \label{H1}
If the second integral cohomology group $H^2_\Z$ and the homogeneous polynomial
of degree $2n-2$ on $H^2_\Z$ defined by the first Pontryagin class are 
given, then there exist at most finitely many diffeomorphism types of 
compact hyperk\"ahler manifolds of complex dimension $2n$ realizing this 
structure. 

\hfill

Equivalently, instead of fixing the first Pontryagin class, one can fix 
the Beauville-Bogomolov-Fujiki form $q$, and thus give the abelian group 
$H^2_\Z$ a ring structure. Also, instead of fixing the whole intersection 
form $q$, we can fix just the following two topological invariants: the 
Fujiki constant $c$ and the discriminant $d$ of the Beauville-Bogomolov-Fujiki 
(BBF) form $q$. The author has established the following finiteness result in 
\cite{_Kamenova_} (for fibrations) and in \cite{_Kamenova2_} (in general).  

\hfill 

\theorem 
Fix a positive rational number $c$ and an integer number $d$. Then 
there are at most finitely many deformation classes of hyperk\"ahler manifolds 
with Fujiki constant $c$ and discriminant $d$ of the BBF lattice $(\Lambda,q)$.

\hfill

For a given diffeomorphic structure on a manifold $M$, Huybrechts proved 
the existence of only finitely many deformation types of hyperk\"ahler 
metrics on $M$ (Theorem 2.1 in \cite{_Huybrechts:finiteness_}). 

\hfill

\theorem (Huybrechts, \cite{_Huybrechts:finiteness_}) \label{H2}
Let $M$ be a fixed compact manifold. Then there exist at most finitely
many different deformation types of irreducible holomorphic symplectic 
complex structures on $M$.

\hfill 

A deformation of a hyperk\"ahler Lagrangian fibration is 
a deformation of the pair $(M, L)$, where $L$ is the line bundle associated 
to the fibration, obtained by pulling back an ample class from the base. 
The author together with Misha Verbitsky had established the following 
finiteness results in \cite{_KV:fibrations_} about Lagrangian fibrations. 

\hfill

\theorem (Kamenova-Verbitsky, \cite{_KV:fibrations_})
Let $M$ be a fixed compact complex manifold of dimension $2n$ and 
$b_2(M) \geq 7$. Then there are 
only finitely many deformation types of hyperk\"ahler 
Lagrangian fibrations on $M$. 

\hfill

The proof of our theorem above relies on Ananin-Verbitsky's general result 
about lattices, which we apply to the lattice $(H^2(M, \Z),q)$, and some 
of its sublattices. 

\hfill

\proposition \cite[Proposition 3.2]{_Ananin_Verbitsky_} 
\label{_lattices_Prop_} 
{\sl Let\/ $V$ be an\/ $\mathbb R$-vector space equipped with a 
non-degene\-rate symmetric form of signature\/ $(s_+,s_-)$ with\/ 
$s_+\ge3$ and\/ $s_-\ge1$. Consider a lattice\/ $L\subset V$. Let\/ 
$\Gamma$ be a subgroup of finite index in\/ $\text{O}(L)$, and\/ $l\in L$. 
Then\/ $\Gamma\cdot\text{Gr}_{++}(l^\perp)$ is dense in $\text{Gr}_{++}(V)$.} 

\hfill

Using Charles's boundedness result for families of hyperk\"ahler varieties 
\cite{_Charles_}, which sharpens Koll\'ar-Matsusaka's classical theorem, 
we were able to prove 
the following finiteness result for Lagrangian fibrations in \cite{_Kamenova_}

\hfill

\theorem 
Consider a Lagrangian fibration $\pi:M \arrow\C P^n$ such that 
there is a line bundle $P$ on $M$ with $q(P)>0$ and with a given $P$-degree 
$d$ on the general fiber $F$ of $\pi$, i.e., $P^n \cdot F = d$. 
Then there are at most finitely many deformation classes of hyperk\"ahler 
manifolds $M$ with a fibration structure as above. 

\hfill

A holomorphic Lagrangian fibration $\pi:M \arrow B$ is stable if the 
discriminant divisor $D$ is a submanifold of the base $B$, 
and each singular fiber is stable (i.e., it is reduced and it is of type 
$I_k$, where $1 \leq k \leq \infty$). Let's assume that the 
base $B$ is smooth, and hence $B = \C P^n$ by Hwang's theorem, 
\cite{_Hwang:base_}. The main result of this paper is the following. 

\hfill

\theorem 
Consider a stable Lagrangian fibration $\pi:M \arrow \C P^n$ with 
discriminant divisor $D \subset \C P^n$. 
Then there are only finitely many Lagrangian fibrations 
$\pi:M' \arrow \C P^n$ with discriminant divisor 
$D \subset \C P^n$.


\section{Definitions and classical results in hyperk\"ahler geometry}


\definition  A {\bf hyperk\"ahler manifold} is a compact complex simply 
connected K\"ahler manifold with a holomorphic symplectic form. 
A hyperk\"ahler manifold $M$ is called {\bf irreducible} if 
$H^{2,0}(M)=\C \sigma$.

\hfill

In this paper we consider irreducible hyperk\"ahler manifolds, 
and for simplicity, we shall refer to them as hyperk\"ahler manifolds. 
According to Matsushita's theorem, there are a lot of restrictions on the 
possible fibration structures on hyperk\"ahler manifolds. 

\hfill

\theorem
(Matsushita, \cite{_Matsushita:fibred_}). \label{Matsushita}
Let $\pi:\; M \rightarrow B$ be a surjective holomorphic map with 
connected fibers 
from a hyperk\"ahler manifold $M$ to a base $B$, with $0<\dim B < \dim M$.
Then $\dim B = \frac{1}{2} \dim M$, and the fibers of $\pi$ are 
holomorphic Lagrangian abelian varieties (i.e., the symplectic 
form vanishes when restricted to the fibers: $\sigma|_F=0$).

\hfill

Such a map $\pi:\; M \rightarrow B$ as above is called a 
{\bf (holomorphic) Lagrangian fibration}.
The set of critical values of $\pi$ is called the {\bf discriminant locus} 
$D \subset B$, and Hwang-Oguiso \cite{HO1} have shown that $D$ is a 
hypersurface if it is non-empty. 

\hfill

\remark 
Matsushita proved that if $B$ is smooth and 
$M$ is projective, then $B$ has the same rational cohomology as $\C P^n$. 
Hwang (\cite{_Hwang:base_}) was able to improve this statement, 
and proved that under the smoothness hypothesis $B\cong \C P^n$. 

\hfill

In \cite{HO1}, Hwang and Oguiso describe the structure of a general singular 
fiber of a holomorphic Lagrangian fibration $\pi:\; M \rightarrow B$, 
similar to Kodaira's classification of singular fibers of elliptic 
fibrations on surfaces. A general singular fiber can be transformed into a 
singular fiber of stable type by the method of stable reduction, \cite{HO2}. 

\hfill 

\definition
A singular fiber of a holomorphic Lagrangian fibration $\pi:\; M \rightarrow B$ 
is of {\bf stable type} if its characteristic cycles are of type $I_k$, 
where $0 \leq k \leq \infty$, and $I_\infty \doteq A_\infty$. We say that $\pi$ 
is a {\bf stable Lagrangian fibration} if the discriminant $D \subset B$ is 
a submanifold and each singular fiber $f^{-1} (b), b \in D$ is reduced and of 
stable type. 

\hfill

For a hyperk\"ahler manifold $M$, on $H^2(M, \Z)$ there is a natural 
primitive integral quadratic form of signature $(3, b_2-3)$, called the 
{\bf  Beauville-Bogomolov-Fujiki form}, or {\bf BBF form}, see 
\cite{_Beauville_} and \cite{_Huybrechts:lec_}. 
One way to define it is using the Fujiki relation.

\hfill

\theorem
(Fujiki, \cite{_Fujiki:HK_}) \label{Fujiki_formula}
Let $\eta\in H^2(M, \Z)$ and $\dim M=2n$, where $M$ is a 
hyperk\"ahler manifold. Then $\eta^{2n}= c \cdot q(\eta,\eta)^n$,
for a primitive integral quadratic form $q$ on $H^2(M, \Z)$, where $c>0$ is a 
constant depending on the topological type of $M$. The constant $c$ in Fujiki's 
formula is called the {\bf Fujiki constant}.

\section{Main result}

In this section we prove the main theorem regarding finiteness of stable 
fibrations with a given discriminant divisor. 

\hfill

\theorem
Consider a stable Lagrangian fibration $\pi:M \arrow \C P^n$ with 
discriminant divisor $D \subset \C P^n$. 
Then there are only finitely many Lagrangian fibrations 
$\pi:M' \arrow \C P^n$ with discriminant divisor 
$D \subset \C P^n$.

\hfill

{\bf Proof}. 
After a degenerate twistor deformation (as in Verbitsky's 
\cite[Theorem 1.10]{_Verbitsky:degenerate_}), without any loss of 
generality we can assume that $M$ is projective, because can always replace 
the original fibration $\pi: M \arrow \C P^n$ 
with another fibration with holomorphically the same fibers, the same 
base, and projective total space by Verbitsky's theorem. Let $L$ be an 
ample line bundle on $M$. By Matsushita's \ref{Matsushita}, 
the smooth fibers $M_b$ over $b \in \C P^n \setminus D$ are abelian 
varieties of dimension $n$. By Oguiso's result \cite{Oguiso}, the 
restriction map on cohomology $H^2(M, \Z) \arrow H^2(M, \Z)$ 
has rank one, and therefore, the fibers $M_b$ are polarized 
by the restriction of the ample line bundle $L$. 
Over the complement of the discriminant divisor, the restriction map 
$M \setminus \pi^{-1} (\C P^n \setminus D) \arrow \C P^n \setminus D$, 
is a holomorphic fibration of polarized abelian varieties of dimension $n$. 
After passing through a level structure, we obtain a holomorphic moduli map  
$m: \C P^n \setminus D \arrow {\cal A}_{n, \alpha}$, where ${\cal A}_{n, \alpha}$ 
is the moduli space of $n$-dimensional abelian varieties with a polarization 
of type $\alpha$, where $\alpha$ is the polarization coming from the given 
ample line bundle $L$. Notice that since $M$ is projective, there are local 
(in the analytic topology) sections from $\C P^n \setminus D$ to $M$, 
which can be glues together by the Cousin problem, because 
$\C P^n \setminus D$ is a Stein variety, and thus the map to 
${\cal A}_{n, \alpha}$ is well defined. 

\hfill

Since we started with a stable Lagrangian fibration, the divisor $D$ is 
smooth, and therefore, one can construct a standard complete 
metric on $\C P^n \setminus D$ with finite volume (see \cite[Section 2]{JY} 
for an explicit construction). For the standard construction of a complete 
finite volume metric on $X \setminus R$, one needs $X$ to be a compact 
K\"ahler manifold, and $R$ to be a divisor on $X$ with at most simple 
normal crossing singularities. In our case, $X = \C P^n$ and $R = D$. 

\hfill

Consider the moduli map  
$m: \C P^n \setminus D \arrow {\cal A}_{n, \alpha}$. 
The orbifold universal cover of ${\cal A}_{n, \alpha}$ is the 
Siegel upper-half space $\mathbb{H}_n$, which carries a natural 
invariant metric as 
$\mathbb{H}_n = \Sp(2n, \Z) \backslash \Sp(2n, \R) / \U(n, \C)$. 
The holomorphic sectional curvature of ${\cal A}_{n, \alpha}$, being a local 
invariant, is the same as the holomorphic sectional curvature of its \'etale 
cover $\mathbb{H}_n$, and thus, it is bounded from above by $-1$. 
By Royden's generalized Schwartz lemma for K\"ahler manifolds \cite{R},  
any such family $m: \C P^n \setminus D \arrow {\cal A}_{n, \alpha}$ is 
equicontinuous. Then, applying Borel's extension theorem \cite{B}, 
any such map $m$ extends to a holomorphic map 
$\C P^n \arrow \bar{\cal A}_{n, \alpha}$ into the Satake compactification 
$\bar{\cal A}_{n, \alpha}$. 
Equicontinuity implies that such extended maps 
$\C P^n \arrow \bar{\cal A}_{n, \alpha}$ can only lie in finitely many homotopy 
classes by \cite[Theorem 2.1 and Section 7]{JY}. On the other hand, any two 
homotopic maps into ${\cal A}_{n, \alpha}$ that extend to maps into 
$\bar{\cal A}_{n, \alpha}$ coincide by \cite[Lemma 3.1 and 
Section 7]{JY} (see also \cite[Step 3 in Theorem 1.1]{Farb}), 
because ${\cal A}_{n, \alpha}$ has strictly negative sectional curvature. 
Therefore, we obtain finiteness of fibrations $\pi:M' \arrow \C P^n$ with a 
given discriminant divisor $D \subset \C P^n$. \endproof

\hfill

\conjecture
We conjecture that one can fix even less data than fixing the whole divisor 
$D \subset \C P^n$, and still obtain finiteness. 
Namely, there should be only finitely many Lagrangian fibrations 
$\pi:M \arrow \C P^n$ with a given degree $d$ of the discriminant divisor 
$D \subset \C P^n$, where the total space $M$ is hyperk\"ahler. 

\hfill

\noindent{\bf Acknowledgments.} 
The author is grateful to Prof. Tian for their past conversations about 
generalizations of Schwartz's lemma, and about Jost-Yau's rigidity results. 
She also thanks Sam Grushevsky for their conversations about moduli spaces 
of Abelian varieties, Ariyan Javanpeykar for his sugestions and improvements, 
and special thanks to the referee for their suggestions. 

\hfill

\noindent{\bf Conflict of interest statement.} 
On behalf of all authors, the corresponding author states that there is no 
conflict of interest.


{\small

\hfill

\noindent {\sc Ljudmila Kamenova\\
Department of Mathematics, 3-115 \\
Stony Brook University \\
Stony Brook, NY 11794-3651, USA,} \\
\tt kamenova@math.stonybrook.edu\\

\end{document}